\documentclass[a4paper,12pt,twoside]{article}
\usepackage{amsmath,amsthm,amssymb}
\usepackage[usenames]{color}
\numberwithin{equation}{section}

\makeatletter

\ifcase\@ptsize
\def\rpkern{\mathchoice{\kern-1.45em}{\kern-1.11em}{}{}}%
\def\grpkern{\mathchoice{\kern-1.013em}{\kern-0.825em}{}{}}%
\or
\def\rpkern{\mathchoice{\kern-1.44em}{\kern-1.11em}{}{}}%
\def\grpkern{\mathchoice{\kern-1.00em}{\kern-0.81em}{}{}}%
\or
\def\rpkern{\mathchoice{\kern-1.472em}{\kern-1.14em}{}{}}%
\def\grpkern{\mathchoice{\kern-1.00em}{\kern-0.815em}{}{}}%
\fi

\def\minibullet{\mathchoice%
{\raise0.2ex\hbox{$\scriptstyle\bullet$}}%
{\raise0.26ex\hbox{$\scriptscriptstyle\bullet$}}{}{}}
\def\butabullet{\mathchoice%
{\raise0.8ex\hbox{$\scriptstyle\bullet$}{\kern-0.365em}%
\lower0.4ex\hbox{$\scriptstyle\bullet$}}%
{\raise0.75ex\hbox{$\scriptscriptstyle\bullet$}{\kern-0.335em}%
\lower0.25ex\hbox{$\scriptscriptstyle\bullet$}}{}{}}

\def\customprod#1#2%
{\operatorname{\coprod\kern#2{#1}\kern#2\prod}\displaylimits}

\newcommand{\bC}{\mathbb{C}}

\newcommand{\bN}{\mathbb{N}}

\newcommand{\bR}{\mathbb{R}}

\newcommand{\bZ}{\mathbb{Z}}

\newcommand{\cA}{\mathcal{A}}

\newcommand{\cM}{\mathcal{M}}

\newcommand{\fS}{\mathfrak{S}}

\renewcommand{\a}{\alpha}

\newcommand{\g}{\gamma}
\renewcommand{\d}{\delta}
\newcommand{\e}{\varepsilon}
\newcommand{\z}{\zeta}

\newcommand{\m}{\mu}

\renewcommand{\r}{\rho}
\newcommand{\s}{\sigma}
\renewcommand{\t}{\tau}

\newcommand{\vp}{\varphi}

\renewcommand{\Re}{\mathrm{Re}\,}

\renewcommand{\det}{\mathrm{det}}

\newcommand{\ord}{\mathrm{ord}}

\newcommand{\Gauss}[1]{\lfloor{#1}\rfloor}

\newcommand{\boldtitle}[1]{\title{\bfseries #1}}

\newcommand{\layout}{%
\topmargin -1in
\oddsidemargin -1in
\evensidemargin -1in
\hoffset -0.07in
\voffset 0.00in
\newlength{\sidemargin}
\setlength{\sidemargin}{-2in}
\addtolength{\sidemargin}{\paperwidth}
\addtolength{\sidemargin}{-\textwidth}
\divide \sidemargin by 2
\setlength{\oddsidemargin}{\sidemargin}
\setlength{\evensidemargin}{\sidemargin}
\setlength{\topmargin}{-2in}
\addtolength{\topmargin}{\paperheight}
\addtolength{\topmargin}{-\textheight}
\addtolength{\topmargin}{-\headheight}
\addtolength{\topmargin}{-\headsep}
\divide \topmargin by 2}

\newenvironment{MSC}{%
\smallbreak
\noindent 2000\ \textbf{Mathematics Subject Classification\,:}\ 
}

\newenvironment{keywords}{%
\noindent\textbf{Key words and phrases\,:}\itshape}

\theoremstyle{theorem}
\newtheorem*{multitheorem}{\variable@name}

\theoremstyle{definition}
\newcommand{\variable@name}{Theorem}
\newtheorem*{multiproclaim}{\variable@name}

\theoremstyle{plain}
\newtheorem{thm}{Theorem}[section]
\newtheorem{prop}[thm]{Proposition}
\newtheorem{lem}[thm]{Lemma}
\newtheorem{cor}[thm]{Corollary}

\theoremstyle{definition}

\newtheorem{example}[thm]{Example}
\newtheorem{remark}[thm]{Remark}

\newcommand{\bsym}{\boldsymbol}

\newcommand{\braket}[1]{\langle{#1}\rangle}

\newcommand{\sgn}{\mathrm{sgn}}
\newcommand{\Gau}[1]{\lfloor{#1}\rfloor}

\renewcommand{\gcd}[1]{\mathrm{gcd}(#1)}
\newcommand{\lcm}[1]{\mathrm{lcm}(#1)}

\renewcommand{\pmod}[1]{(\mathrm{mod}\ {#1})}
\newcommand{\wt}[1]{\widetilde{#1}}

\boldtitle{Arithmetical properties of Multiple Ramanujan sums} 
\author{Yoshinori YAMASAKI}
\date{\today}

\layout

\begin{document}

\setlength{\baselineskip}{14.9pt}
\maketitle

\begin{abstract}
 In the present paper, we introduce a multiple Ramanujan sum for
 arithmetic functions, which gives a multivariable extension of the
 generalized Ramanujan sum studied by D. R. Anderson and
 T. M. Apostol. We then find fundamental arithmetic properties of the
 multiple Ramanujan sum and study several types of Dirichlet series
 involving the multiple Ramanujan sum. As an application, we
 evaluate higher-dimensional determinants of higher-dimensional
 matrices, the entries of which are given by values of the multiple
 Ramanujan sum.
\begin{MSC}
 {\it Primary} 11A25; {\it Secondary} 11C20.
\end{MSC} 
\end{abstract}

\section{Introduction}
\label{sec:Introduction}

 In 1918, S. Ramanujan \cite{Ramanujan1918} studied the sum  
\begin{equation}
\label{def:Ramanujan sum} 
 c(k,n):
=\sum_{l\,\pmod{k} \atop \gcd{l,k}=1}e(k,nl)
=\sum_{d|\,\gcd{k,n}}\m\bigl(\frac{k}{d}\bigr)d, 
\end{equation}
 where, in the first sum, $l$ runs over a reduced residue system modulo
 $k$ with $\gcd{l,k}=1$, $e(r,n):=\exp(2\pi\sqrt{-1}n/r)$ and $\m$ is
 the M\"obius function. The sum $c(k,n)$ is called the Ramanujan sum
 (or the Ramanujan trigonometric sum) and is widely investigated in
 connection with, for example, even arithmetic functions
 \cite{Cohen1955,Cohen1960} and cyclotomic polynomials
 \cite{Motose2005,Nicol1962}. See \cite{McCarthy1986b} for
 the arithmetic theory of the Ramanujan sum. Among several 
 generalizations and variations of $c(k,n)$, D. R. Anderson and
 T. M. Apostol \cite{AndersonApostol1953} (see also \cite{Apostol1972}) 
 considered the sum
\begin{equation}
\label{def:generalizedRS}
 S_{f,g}(k,n):=\sum_{d|\,\gcd{k,n}}f\bigl(\frac{k}{d}\bigr)g(d),
\end{equation}
 where $f$ and $g$ are arithmetic functions. Clearly, $S_{f,g}(k,n)$
 extends the right-most expression in formula \eqref{def:Ramanujan sum}
 and hence gives a generalization of the Ramanujan sum.

 Motivated by the study of the above generalized Ramanujan sum, in the
 present paper, we examine the following type of multiple sum for
 arithmetic functions $f_1,\ldots,f_{m+1}$;
\[
 S_{f_1,\ldots,f_{m+1}}(n_1,\ldots,n_{m+1})
:=\sum_{d_j|\,\gcd{n_1,\ldots,n_{j+1}}\atop (j=1,\ldots,m)}
 f_1\bigl(\frac{n_1}{d_1}\bigr)f_2\bigl(\frac{d_1}{d_2}\bigr)\cdots
 f_m\bigl(\frac{d_{m-1}}{d_{m}}\bigr)f_{m+1}\bigl(d_m\bigr),
\]
 where $\gcd{n_1,\ldots,n_{j+1}}$ is the greatest common divisor of
 $n_1,\ldots,n_{j+1}$. We call this a {\it multiple Ramanujan sum} for 
 $f_1,\ldots,f_{m+1}$. Notice that the above expression gives the
 generalized Ramanujan sum \eqref{def:generalizedRS} when $m=1$ and,
 moreover, the Dirichlet convolution $f_1*\cdots *f_{m+1}$ of
 $f_1,\ldots,f_{m+1}$ in the ``diagonal case''; $n_1=\cdots=n_{m+1}$. 

 The present paper is organized as follows. In
 Section~\ref{sec:definition}, we introduce a multiple Ramanujan sum
 $S^{(\g_1,\ldots,\g_m)}_{f_1,\ldots,f_{m+1}}$ with positive integer
 parameters $\g_1,\ldots,\g_m$ so that
 $S_{f_1,\ldots,f_{m+1}}=S^{(1,\ldots,1)}_{f_1,\ldots,f_{m+1}}$ and
 study its fundamental properties as a multivariable arithmetic
 function, such as the degeneracies and the multiplicativity (see
 \cite{Vaidyanathaswamy1931} for the theory of multivariable arithmetic
 functions). Then, since $S^{(\g_1,\ldots,\g_m)}_{f_1,\ldots,f_{m+1}}$
 belongs to the class of even arithmetic functions $\pmod{n_1}$ as a
 function of $n_2,\ldots,n_{m+1}$ in the sense of Cohen
 \cite{Cohen1960}, we calculate its finite Fourier expansion, which any
 even arithmetic function possesses. This expression is important with
 respect to Section~\ref{sec:determinant}.
 In Section~\ref{sec:Dirichlet series}, we study several types of
 Dirichlet series having coefficients that are given by
 $S^{(\g_1,\ldots,\g_m)}_{f_1,\ldots,f_{m+1}}$. We treat not only 
 single-variable Dirichlet series, but also multivariable Dirichlet
 series. For instance, an analogue of the formula of J. M. Borwein and
 K. K. Choi \cite{BorweinChoi2003}, which contains the classical
 Ramanujan formula concerning the divisor function $\s_{a}$, is
 obtained. Section~\ref{sec:determinant} is devoted to a
 higher-dimensional generalization of the so-called Smith determinant
 \cite{Smith187576}. We evaluate higher-dimensional determinants 
 (hyperdeterminants) of higher-dimensional matrices (hypermatrices),
 the entries of which are given by values of
 $S^{(\g_1,\ldots,\g_m)}_{f_1,\ldots,f_{m+1}}$. In fact, we derive a
 hyperdeterminant formula for even multivariable arithmetic
 functions. This includes the results of P. J. McCarthy
 \cite{McCarthy1986a} and K. Bourque and S. Ligh \cite{BourqueLigh1993}
 for the $2$-dimensional case, that is, the usual determinant case, and
 partially of P. Haukkanen \cite{Haukkanen1992}, for the
 higher-dimensional case. 
 
 We use the following notations in the present paper. The set of natural 
 numbers, the ring of rational integers, the field of real numbers and
 the field of complex numbers are denoted respectively as $\bN$, $\bZ$,
 $\bR$, and $\bC$. For $n_1,\ldots,n_{k}\in\bN$, $\gcd{n_1,\ldots,n_k}$ 
 (resp. $\lcm{n_1,\ldots,n_k}$) represents for the greatest common
 divisor (resp. least common multiple) of $n_1,\ldots,n_k$. For
 $x\in\bR$, $\Gau{x}$ is the greatest integer not exceeding $x$. We
 denote the M\"obius function as $\m(n)$, the Euler totient function as
 $\vp(n)$, the power function as $\d^{x}(n):=n^{x}$ for $x\in\bC$, and
 the identity element in the ring of arithmetic functions with respect to
 the Dirichlet convolution $*$ as $\e(n):=\Gau{\frac{1}{n}}=1$ if $n=1$,
 and $0$ otherwise. Note that $\d^{0}*\m=\e$. Throughout the present
 paper, we consider a product (resp. a sum) over an empty set to always
 be equal to $1$ (resp. $0$).

\section{A multiple Ramanujan sum}
\label{sec:definition}

\subsection{Preliminary:\ $\bsym{\g}$-convolutions}

 Let $\cA$ be the set of all complex-valued arithmetic functions
 $f:\bN\to\bC$. We always understand that $f(x)=0$ if $x\notin\bN$ for
 $f\in\cA$. The set of all multiplicative and completely multiplicative
 arithmetic functions are respectively denoted by $\cM$ and
 $\cM^{c}$. Namely, $f\in\cM$ (resp. $f\in\cM^c$) means that
 $f(mn)=f(m)f(n)$ for $\gcd{m,n}=1$ (resp. for all $m,n\in\bN$). As
 usual, the product $fg\in\cA$ of $f,g\in\cA$ is defined by
 $(fg)(n):=f(n)g(n)$.
 
 Let $\g\in\bN$. We define the $\g$-convolution of $f,g\in\cA$ by
\[
 (g*_{\g}f)(n)
:=\sum_{d^{\g}|n}f\bigl(\frac{n}{d^{\g}}\bigr)g(d^{\g}).
\]
 In particular, $*_{1}=*$ stands for the usual Dirichlet convolution.
 Define the function $a_{\g}\in\cA$ by $a_{\g}(n)=1$ if there exists
 $d\in\bN$ such that $n=d^{\g}$, and $0$ otherwise. Then, it is clear
 that $g*_{\g}f=g^{[\g]}*f$, where $g^{[\g]}:=a_{\g}g$. Note that
 the product $*_{\g}$ does not satisfy the commutativity or 
 associativity properties unless $\g=1$. We therefore inductively define
 the $\bsym{\g}=(\g_1,\ldots,\g_m)$-convolution for
 $\g_1,\ldots,\g_m\in\bN$ of $f_1,\ldots,f_{m+1}\in\cA$ as
\begin{align*}
 (f_{m+1}*_{\g_{m}}\cdots *_{\g_1}f_{1})(n)
:&=\bigl((f_{m+1}*_{\g_{m}}\cdots *_{\g_{2}}f_{2})*_{\g_1}f_{1}\bigr)(n)\\
&=\sum_{d_{m}^{\g_m}|\cdots|d_{1}^{\g_{1}}|n}
f_{1}\bigl(\frac{n}{d_1^{\g_1}}\bigr)f_{2}\bigl(\frac{d_1^{\g_1}}{d_{2}^{\g_2}}\bigr)\cdots f_{m}\bigl(\frac{d_{m-1}^{\g_{m-1}}}{d_m^{\g_m}}\bigr)f_{m+1}\bigl(d_m^{\g_m}\bigr).
\end{align*}

 Let $L(s;f):=\sum^{\infty}_{n=1}f(n)n^{-s}$ be the Dirichlet series
 attached to $f\in\cA$ and let $\s(f)\in\bR\cup\{\infty\}$ be the
 abscissa of absolute convergence of $L(s;f)$. Define
 $f^{\braket{\g}}\in\cA$ for $\g\in\bN$ by
 $f^{\braket{\g}}(n):=f(n^{\g})$. Then, it is clear that
 $L(\g s;f^{\braket{\g}})=L(s;f^{[\g]})$ with
 $\s(f^{[\g]})\le\s(f)$. For $f_1,\ldots,f_{m+1}\in\cA$, we here
 calculate the Dirichlet series attached to
 $f_{m+1}*_{\g_{m}}\cdots *_{\g_1}f_{1}$. 

\begin{prop}
\label{prop:Dirichlet}
 Suppose $\g_{0}|\g_{1}|\g_{2}|\cdots|\g_{m}$ with $\g_0=1$. Then, we
 have for $\Re(s)>\max\{\s(f_{1}),\ldots,\s(f_{m+1})\}$
\begin{equation}
\label{for:Dirichletseries}
 L\bigl(s;f_{m+1}*_{\g_{m}}\cdots *_{\g_{1}}f_{1}\bigr)
=\prod^{m+1}_{j=1}L\bigl(s;f^{[\g_{j-1}]}_{j}\bigr).
\end{equation}
\end{prop}
\begin{proof}
 Assume $d^{\g_{m}}_{m}|\cdots|d^{\g_{1}}_{1}|n$. Then, since
 $\g_{1}|\g_{2}|\cdots|\g_{m}$, we can write
 $d^{\g_{j}}_{j}=d^{\g_{j+1}}_{j+1}k^{\g_{j}}_{j+1}$ for $1\le j\le m-1$
 and
 $n=k_1d_1^{\g_1}=\cdots=k_{1}k^{\g_{1}}_{2}\cdots k^{\g_{m-1}}_{m}d^{\g_m}_m$ 
 with $k_{1},k_{2},\ldots,k_{m}\in\bN$. We then write $d_{m}=k_{m+1}$,
 and so
\begin{align*}
 L\bigl(s;f_{m+1}*_{\g_{m}}\cdots *_{\g_1}f_{1}\bigr)
&=\sum^{\infty}_{k_1,\ldots,k_{m+1}=1}f_1(k_1)f_2(k_2^{\g_1})\cdots
 f_{m+1}(k_{m+1}^{\g_{m}})\bigl(k_1k_2^{\g_1}\cdots k_{m+1}^{\g_m}\bigr)^{-s}\\
&=L\bigl(s;f_1\bigr)L\bigl(\g_1 s;f_2^{\braket{\g_2}}\bigr)\cdots L\bigl(\g_m
 s;f_{m+1}^{\braket{\g_m}}\bigr).
\end{align*}
 This ends the proof.
\end{proof}

\begin{remark}
 We cannot expect that $L(s;f_{m+1}*_{\g_{m}}\cdots *_{\g_{1}}f_{1})$
 is expressed as a product of a Dirichlet series such as 
 \eqref{for:Dirichletseries} for general
 $\bsym{\g}=(\g_1,\ldots,\g_m)\in\bN^{m}$.
\end{remark}

 In the next subsection, we will introduce a multiple Ramanujan sum for 
 $f_1,\ldots,f_{m+1}\in\cA$, which gives the $\bsym{\g}$-convolution
 $f_{m+1}*_{\g_{m}}\cdots *_{\g_{1}}f_{1}$ in the diagonal case (for
 more detail, see Proposition~\ref{prop:elementary} $(\mathrm{iv})$).

\subsection{Definition of $S^{\bsym{\g}}_{\bsym{f}}$ and its basic properties} 

 Let $\bsym{\g}=(\g_1,\ldots,\g_m)\in\bN^{m}$ and
 $\bsym{f}=(f_1,\ldots,f_{m+1})\in\cA^{m+1}$. We define a multiple
 Ramanujan sum 
 $S^{\bsym{\g}}_{\bsym{f}}=S^{(\g_1,\ldots,\g_m)}_{f_1,\ldots,f_{m+1}}$
 of $\bsym{f}$ with the parameter $\bsym{\g}$ by 
\[
 S^{\bsym{\g}}_{\bsym{f}}(n_1,\ldots,n_{m+1})\\
:=\sum_{d^{\g_j}_j|\,\gcd{n_1,\ldots,n_{j+1}} \atop (j=1,\ldots,m)}
f_{1}\bigl(\frac{n_1}{d_1^{\g_1}}\bigr)f_{2}\bigl(\frac{d_1^{\g_1}}{d_{2}^{\g_2}}\bigr)\cdots f_{m}\bigl(\frac{d_{m-1}^{\g_{m-1}}}{d_m^{\g_m}}\bigr)f_{m+1}\bigl(d_m^{\g_m}\bigr),
\]
 where the sum is taken over all $m$-tuples $(d_1,\ldots,d_m)\in\bN^{m}$ 
 satisfying $d^{\g_j}_j|\,\gcd{n_1,\ldots,n_{j+1}}$ for each
 $1\le j\le m$. Note that the summand vanishes unless
 $d_m^{\g_m}|\cdots |d_1^{\g_1}|n_1$. We write
 $S_{\bsym{f}}:=S^{(1,\ldots,1)}_{\bsym{f}}$ and understand that
 $S_{f}(n)=f(n)$ when $m=0$. These are the fundamental properties of
 $S^{\bsym{\g}}_{\bsym{f}}$, which are obtained from the elementary
 properties of the gcd-function.

\begin{prop}
\label{prop:elementary}
 $(\mathrm{i})$\ For any $1\le j\le m+1$, we have 
\begin{equation}
\label{for:Srec}
 S^{\bsym{\g}}_{\bsym{f}}(n_1,\ldots,n_{m+1})
=S^{(\g_1,\ldots,\g_{j-1})}_{f_1,\ldots,f_{j-1},S^{(\g_j,\ldots,\g_m)}_{f_j,\ldots,f_{m+1}}(\,\cdot\,,n_{j+1},\ldots,n_{m+1})}(n_1,\ldots,n_{j}).
\end{equation}

 $(\mathrm{ii})$\ For any $1\le j\le m+1$, we have
\begin{equation*}
\label{for:degeneracy1}
 S^{\bsym{\g}}_{\bsym{f}}(n_1,\ldots,n_{j-1},1,n_{j+1},\ldots,n_{m+1})
=f_{j}(1)\cdots f_{m+1}(1)\cdot
 S^{(\g_1,\ldots,\g_{j-2})}_{f_1,\ldots,f_{j-1}}(n_1,\ldots,n_{j-1}).
\end{equation*}
 In particular, we have
 $S^{\bsym{\g}}_{\bsym{f}}(1,n_2,\ldots,n_{m+1})=f_{1}(1)\cdots f_{m+1}(1)$. 

 $(\mathrm{iii})$\ For any $1\le j\le m$, we have
\begin{multline}
\label{for:degeneracy2}
 S^{(\g_1,\ldots,\g_{m})}_{f_1,\ldots,f_{j-1},\e,f_{j+1},\ldots,f_{m+1}}(n_1,\ldots,n_{m+1})\\
= S^{(\g_1,\ldots,\g_{j-1},\g_{j+1},\ldots,\g_{m})}_{f_1,\ldots,f_{j-1},f_{j+1},\ldots,f_{m+1}}(n_1,\ldots,n_{j-1},\gcd{n_j,n_{j+1}},n_{j+2},\ldots,n_{m+1})
\end{multline}
 and $S^{(\g_1,\ldots,\g_{m})}_{f_1,\ldots,f_m,\e}(n_1,\ldots,n_{m+1})=S^{(\g_1,\ldots,\g_{m-1})}_{f_1,\ldots,f_m}(n_1,\ldots,n_{m})$.

 $(\mathrm{iv})$\ Let $n_1|n_j$ for all $2\le j\le m+1$. Then, we have 
\begin{equation}
\label{for:genecon}
 S^{\bsym{\g}}_{\bsym{f}}(n_1,\ldots,n_{m+1})
=(f_{m+1}*_{\g_{m}}\cdots *_{\g_1}f_1)(n_1).
\end{equation}\qed 
\end{prop}

\begin{example}
\label{ex:divisor} 
 Let $\bsym{a}=(a_1,\ldots,a_m)\in\bC^{m}$. Define a multiple divisor
 function $\s^{\bsym{\g}}_{\bsym{a}}$ with the parameter $\bsym{\g}$ by 
\begin{equation}
 \s^{\bsym{\g}}_{\bsym{a}}(n_1,\ldots,n_{m+1})
:=\sum_{
\begin{subarray}{c}
 d_m^{\g_m}|\cdots |d_1^{\g_1}|n_1\\
 d^{\g_j}_j|\,\gcd{n_1,\ldots,n_{j+1}}\ (j=1,\ldots,m) 
\end{subarray}
}
d_1^{\g_1 a_1}\cdots d_{m}^{\g_m a_m}.
\end{equation}
 This is a generalization of the usual divisor function
 $\s_a(n):=\sum_{d|n}d^{a}$; $\s_a(n)=\s^{(1)}_{a}(n,n)$. Then, it holds
 that $\s^{\bsym{\g}}_{\bsym{a}}=S^{\bsym{\g}}_{\bsym{f}}$ with
 $f_j=\d^{a_0+a_1+\cdots+a_{j-1}}$ for $1\le j\le m+1$ where
 $a_0=0$. Similarly, for $\bsym{b}=(b_1,\ldots,b_{m+1})\in\bC^{m+1}$,
 one can see that $\s^{\bsym{\g}}_{\widetilde{\bsym{b}}}(n_1,\ldots,n_{m+1})=n_1^{-b_1}S^{\bsym{\g}}_{\d^{b_1},\ldots,\d^{b_{m+1}}}(n_1,\ldots,n_{m+1})$
 where $\wt{\bsym{b}}:=(b_2-b_1,b_3-b_2,\ldots,b_{m+1}-b_{m})\in\bC^{m}$. 
 Note that the sum  
 $Z^{\bsym{\g}}_{n}(a_1,\ldots,a_m):=\s^{\bsym{\g}}_{\bsym{a}}(n,\ldots,n)$
 is studied in \cite{KimotoKurokawaMatsumotoWakayama2005} and called the
 multiple finite Riemann zeta function. From formula
 \eqref{for:Dirichletseries}, if $\g_{0}|\g_{1}|\g_{2}|\cdots|\g_{m}$
 with $\g_0=1$, we have 
 $L(s;Z^{\bsym{\g}}_{\cdot}(a_1,\ldots,a_m))=\prod^{m+1}_{j=1}\z(\g_{j-1}(s-a_0-\cdots-a_{j-1}))$
 where $\z(s)=L(s;\d^{0})$ is the Riemann zeta function. 
\end{example}

\begin{example}
\label{ex:gcd}
 Let $f\in\cA$. Then, the composition function $f\circ \textrm{gcd}$ of
 $f$ and the gcd-function can be expressed in terms of the multiple
 Ramanujan sum. Actually, from the degeneracy formula
 \eqref{for:degeneracy2}, we have
 $(f\circ \textrm{gcd})(n_1,\ldots,n_{m+1})=S_{\bsym{f}}(n_1,\ldots,n_{m+1})$
 with $f_1=\d^{0}$, $f_2=\cdots=f_m=\e$ and $f_{m+1}=f*\m$.
\end{example}

 We next show the multiplicative property of the multiple Ramanujan
 sum $S^{\bsym{\g}}_{\bsym{f}}$. Recall that an arithmetic function
 $F(n_1,\ldots,n_k)$ of $k$-variables is called multiplicative if
\[
 F(m_{1}n_{1},\ldots,m_{k}n_{k})
=F(m_{1},\ldots,m_{k})\cdot F(n_{1},\ldots,n_{k})
\]
 for relatively prime $k$-tuples $(m_{1},\ldots,m_{k})\in\bN^{k}$ and
 $(n_{1},\ldots,n_{k})\in\bN^{k}$. Here, we say that
 $(m_{1},\ldots,m_{k})$ and $(n_{1},\ldots,n_{k})$ are relatively prime 
 if $\gcd{m_{i},n_{j}}=1$ for all $\le i,j\le k$, and equivalently,
 $\gcd{n_1\cdots n_{k},m_{1}\cdots m_{k}}=1$ (see
 \cite{Vaidyanathaswamy1931}). In this case, we have the following Euler 
 product expression 
\[
 F(n_{1},\ldots,n_{k})
=\prod_{p}F(p^{\a_{p,1}},\ldots,p^{\a_{p,k}}),
\]
 where $\a_{p,j}=\ord_{p}n_j$ for $1\le j\le k$. Note that this is a
 finite product since $F(1,\ldots,1)=1$. 

\begin{prop}
\label{prop:mult}
 The function $S^{\bsym{\g}}_{\bsym{f}}$ is multiplicative if
 $f_1,\ldots,f_{m+1}\in\cM$. 
\end{prop}
\begin{proof}
 Assume $(n_{1},\ldots,n_{m+1})$ and $(k_{1},\ldots,k_{m+1})$ are
 relatively prime. Then, for all $1\le j\le m+1$, 
 $\gcd{n_1k_1,\ldots,n_{j}k_{j}}=\gcd{n_1,\ldots,n_{j}}\cdot \gcd{k_1,\ldots,k_{j}}$
 and $\gcd{n_1,\ldots,n_{j}}$ and $\gcd{k_1,\ldots,k_{j}}$ are
 relatively prime. Hence, the first assertion follows in the same
 manner as the that in the proof of Theorem~$1$ in
 \cite{AndersonApostol1953}. 
\end{proof}

\subsection{Finite Fourier expansions}

 An arithmetic function $F(r;n_1,\ldots,n_k)$ of $k$-variables
 $n_1,\ldots,n_k$ is called periodic $\pmod{r}$ if
 $F(r;n_{1},\ldots,n_{k})=F(r;n'_{1},\ldots,n'_{k})$ whenever
 $n_j\equiv n'_j$ $\pmod{r}$ for all $1\le j\le k$ (see
 \cite{McCarthy1986b} for the case of $k=1$). It is well-known that $F$ is
 periodic $\pmod{r}$ if and only if it has an expression of the form of 
\begin{equation}
\label{for:ffeP} 
 F(r;n_{1},\ldots,n_{k})
=\sum^{r}_{l_1,\ldots,l_k=1}a_{r}(l_1,\ldots,l_k)e(r,n_1l_1)\cdots
e(r,n_kl_k)
\end{equation}
 and the coefficients $a_{r}(l_1,\ldots,l_k)$ are uniquely determined by  
\begin{equation}
\label{for:ffcP} 
 a_{r}(l_1,\ldots,l_k)
=\frac{1}{r^k}\sum^{r}_{n_1,\ldots,n_k=1}F(r;n_1,\ldots,n_{k})e(r,-n_1l_1)\cdots e(r,-n_kl_k).
\end{equation}
 Moreover, $F$ is called even $\pmod{r}$ if
 $F(r;n_{1},\ldots,n_{k})=F(r;\gcd{n_1,r},\ldots,\gcd{n_k,r})$.
 Note that $F$ is periodic $\pmod{r}$ if it is even $\pmod{r}$. Then, as 
 shown by E. Cohen \cite{Cohen1960}, $F$ is even $\pmod{r}$ if and only
 if it has an expression of the form of   
\begin{equation}
\label{for:ffeE} 
 F(r;n_{1},\ldots,n_{k})
=\sum_{d_1,\ldots,d_k|r}\a_{r}(d_1,\ldots,d_k)c(d_1,n_1)\cdots c(d_k,n_k)
\end{equation}
 with
\begin{equation}
\label{for:ffcE} 
 \a_{r}(d_1,\ldots,d_k)
=\frac{1}{r^k}\sum_{\d_1,\ldots,\d_k|r}F(r;\d_1,\ldots,\d_k)c\bigl(\frac{r}{\d_1},\frac{r}{d_1}\bigr)\cdots c\bigl(\frac{r}{\d_k},\frac{r}{d_k}\bigr).
\end{equation}
 We call expressions \eqref{for:ffeP} and \eqref{for:ffeE} finite
 Fourier expansions and coefficients $a_{r}$ and $\a_{r}$ finite
 Fourier coefficients of $F$.  

 By definition, the multiple Ramanujan sum
 $S^{\bsym{\g}}_{\bsym{f}}(n_1,\ldots,n_{m+1})$ is even
 $\pmod{n_1}$ as a function of $m$-variables $n_2,\ldots,n_{m+1}$ for
 any $n_1\in\bN$. Then, let us calculate the finite Fourier expansions
 of $S^{\bsym{\g}}_{\bsym{f}}$. To do so, we add one
 parameter (or, a weight) to $S^{\bsym{\g}}_{\bsym{f}}$. Let $\xi$ be an
 arithmetic function of $m$-variables. Set 
\begin{multline*}
 S^{\bsym{\g},\xi}_{\bsym{f}}(n_1,\ldots,n_{m+1})
=S^{(\g_1,\ldots,\g_m),\xi}_{f_1,\ldots,f_{m+1}}(n_1,\ldots,n_{m+1})\\
:=\sum_{d^{\g_j}_j|\,\gcd{n_1,\ldots,n_{j+1}}\atop
 (j=1,\ldots,m)}\xi(d_1^{\g_1},\ldots,d_{m}^{\g_m})f_{1}\bigl(\frac{n_1}{d_1^{\g_1}}\bigr)f_{2}\bigl(\frac{d_1^{\g_1}}{d_{2}^{\g_2}}\bigr)\cdots f_{m}\bigl(\frac{d_{m-1}^{\g_{m-1}}}{d_m^{\g_m}}\bigr)f_{m+1}\bigl(d_m^{\g_m}\bigr).
\end{multline*}
 Here, $S^{\bsym{\g},\xi}_{\bsym{f}}(n_1,\ldots,n_{m+1})$ is
 again even $\pmod{n_1}$ as a function of $n_2,\ldots,n_{m+1}$ for any
 $n_1\in\bN$ and $S^{\bsym{\g}}_{\bsym{f}}=S^{\bsym{\g},1_m}_{\bsym{f}}$
 where $1_m(d_1,\ldots,d_{m})\equiv 1$. We also write
 $S^{\xi}_{\bsym{f}}:=S^{(1,\ldots,1),\xi}_{\bsym{f}}$. As discussed in
 Proposition~\ref{prop:mult}, 
 $S^{\bsym{\g},\xi}_{\bsym{f}}$ can also be shown to be multiplicative if
 $f_1,\ldots,f_{m+1}$ and $\xi$ are multiplicative.

 For $\bsym{f}=(f_1,\ldots,f_{m+1})$, we set
 ${}^{t}\wt{\bsym{f}}=(\d^{0}f_{m+1},\d^{1}f_{m},\ldots,\d^{m}f_1)$.
 Further, for $n\in\bN$, we set
\[
 {}^{t}\xi^{\bsym{\g}}_{n}(d_1,\ldots,d_m):
=\Bigl(\prod^{m}_{j=1}a_{\g_j}\bigl(\frac{n}{d_{m+1-j}}\bigr)\Bigr)\cdot
\xi\bigl(\frac{n}{d_m},\ldots,\frac{n}{d_1}\bigr).
\]
 Then, we obtain the following theorem. Roughly speaking, the finite
 Fourier coefficients of $S^{\bsym{\g},\xi}_{\bsym{f}}$ can again be written
 as a multiple Ramanujan sum.

\begin{thm}
\label{thm:ffeS}
 Write the finite Fourier expansions of $S^{\bsym{\g},\xi}_{\bsym{f}}$ as 
\begin{align*}
 S^{\bsym{\g},\xi}_{\bsym{f}}(n_1,\ldots,n_{m+1})
&=\sum^{n_1}_{l_2,\ldots,l_{m+1}=1}a^{\bsym{\g},\xi}_{\bsym{f},n_1}(l_2,\ldots,l_{m+1})
e(n_1,n_2l_2)\cdots e(n_1,n_{m+1}l_{m+1})\\
&=\sum_{d_2,\ldots,d_{m+1}|n_1}\a^{\bsym{\g},\xi}_{\bsym{f},n_1}(d_2,\ldots,d_{m+1})
c(d_2,n_2)\cdots e(d_{m+1},n_{m+1}).
\end{align*}
 Then, the finite Fourier coefficients
 $a^{\bsym{\g},\xi}_{\bsym{f},n_1}$ and
 $\a^{\bsym{\g},\xi}_{\bsym{f},n_1}$ of $S^{\bsym{\g},\xi}_{\bsym{f}}$
 are respectively given as
\begin{align}
\label{for:ffcPS}
 a^{\bsym{\g},\xi}_{\bsym{f},n_1}(l_2,\ldots,l_{m+1})
&=n_1^{-m}S^{{}^{t}\xi^{\bsym{\g}}_{n_1}}_{{}^{t}\wt{\bsym{f}}}(n_1,l_{m+1},\ldots,l_{2}),\\
\label{for:ffcES}
 \a^{\bsym{\g},\xi}_{\bsym{f},n_1}(d_2,\ldots,d_{m+1})
&=n_1^{-m}S^{{}^{t}\xi^{\bsym{\g}}_{n_1}}_{{}^{t}\wt{\bsym{f}}}\bigl(n_1,\frac{n_1}{d_{m+1}},\ldots,\frac{n_1}{d_{2}}\bigr).
\end{align}
\end{thm}
\begin{proof}
 We only verify formula \eqref{for:ffcES} (formula \eqref{for:ffcPS} can 
 be similarly obtained). From formula \eqref{for:ffcE},
 $\a^{\bsym{\g},\xi}_{\bsym{f},n_1}(d_2,\ldots,d_{m+1})$ is given as 
\[
 n_1^{-m}\sum_{\d_2,\ldots,\d_{m+1}|n_1}S^{\bsym{\g},\xi}_{\bsym{f}}(n_1,\d_2,\ldots,\d_{m+1})c\bigl(\frac{n_1}{\d_{2}},\frac{n_1}{d_2}\bigr)\cdots c\bigl(\frac{n_1}{\d_{m+1}},\frac{n_1}{d_{m+1}}\bigr).
\]
 Further, by changing the order of the summation, this expression is
 equivalent to   
\begin{multline}
\label{for:mid1}
 n_1^{-m}\sum_{e_m^{\g_m}|\cdots|e_1^{\g_1}|n_1}\xi\bigl(e_1^{\g_1},\ldots,e_{m}^{\g_m}\bigr)f_1\bigl(\frac{n_1}{e_1^{\g_1}}\bigr)f_2\bigl(\frac{e_{1}^{\g_{1}}}{e_2^{\g_2}}\bigr)\cdots
 f_m\bigl(\frac{e_{m-1}^{\g_{m-1}}}{e_m^{\g_m}}\bigr)f_{m+1}\bigl(e_m^{\g_m}\bigr)\\
 \times
 \Bigl(\sum_{e^{\g_1}_1|\d_2|n_1}c\bigl(\frac{n_1}{\d_2},\frac{n_1}{d_2}\bigr)\Bigr)\cdots \Bigl(\sum_{e^{\g_m}_m|\d_{m+1}|n_1}c\bigl(\frac{n_1}{\d_{m+1}},\frac{n_1}{d_{m+1}}\bigr)\Bigr).
\end{multline}
 Here, we use the following identity. Let $e|n$ and $d|n$. Then, we have   
\begin{equation}
\label{for:formula for c}
 \sum_{e|\d|n}c\bigl(\frac{n}{\d},\frac{n}{d}\bigr)=
\begin{cases}
 \displaystyle{\frac{n}{e}} & \textrm{if $d|e$},\\
 0 & \textrm{otherwise.}
\end{cases}
\end{equation}
 Actually, one can obtain this formula from the right-most expression
 \eqref{def:Ramanujan sum} of the Ramanujan sum $c(k,n)$ and formula
 $\d^{0}*\m=\e$. Then, applying formula \eqref{for:formula for c},
 we see that \eqref{for:mid1} can be written as
\begin{multline*}
\quad n_1^{-m}\sum_{d_2|e_1^{\g_1}|n_1}
\sum_{d_3|e_2^{\g_2}|e_1^{\g_1}}\cdots
\sum_{d_{m+1}|e_m^{\g_m}|e_{m-1}^{\g_{m-1}}}
\xi\bigl(e_1^{\g_1},\ldots,e_{m}^{\g_m}\bigr)\\
\times f_1\bigl(\frac{n_1}{e_1^{\g_1}}\bigr)f_2\bigl(\frac{e_{1}^{\g_{1}}}{e_2^{\g_2}}\bigr)\cdots
 f_m\bigl(\frac{e_{m-1}^{\g_{m-1}}}{e_m^{\g_m}}\bigr)f_{m+1}\bigl(e_m^{\g_m}\bigr)\frac{n_1}{e_1^{\g_1}}\cdots \frac{n_1}{e_m^{\g_m}}.\quad
\end{multline*}   
 Changing variables $e'_k=n_1/e_{k}^{\g_k}$ for all $1\le k\le m$, we
 have $a_{\g_k}(n_1/e'_k)=1$. Moreover, it holds that
 $e'_k|\frac{n_1}{d_{k+1}}$ for $1\le k\le m$ and
 $e_1'|\cdots|e'_m|n_1$, because
 $e_m^{\g_m}|\cdots|e_1^{\g_1}|n_1$. Hence, we can rewrite the above
 expression as
\begin{multline}
\label{for:mid2}
 n_1^{-m}\sum_{e'_m|n_1\atop e'_m|\frac{n_1}{d_{m+1}}}
\sum_{e'_{m-1}|e'_m\atop e'_{m-1}|\frac{n_1}{d_{m}}}\cdots
\sum_{e'_1|e'_2\atop e'_1|\frac{n_1}{d_{2}}}
 a_{\g_1}\bigl(\frac{n_1}{e'_1}\bigr)\cdots
 a_{\g_m}\bigl(\frac{n_1}{e'_m}\bigr)
\cdot 
 \xi\bigl(\frac{n_1}{e'_1},\ldots,\frac{n_1}{e'_m}\bigr)\\
 \times f_1\bigl(e'_1\bigr)f_2\bigl(\frac{e'_2}{e'_1}\bigr)\cdots
 f_m\bigl(\frac{e'_m}{e'_{m-1}}\bigr)f_{m+1}\bigl(\frac{n_1}{e'_m}\bigr)e'_1\cdots e'_m.\quad
\end{multline}
 It is easy to see that \eqref{for:mid2} coincides with the right-hand
 side of formula \eqref{for:ffcES}. This completes the proof of theorem.
\end{proof}

\begin{remark}
 Suppose that $\g_0|\g_1|\cdots |\g_m$ with $\g_0=1$. Then, we have   
 $S^{{}^{t}\xi^{\bsym{\g}}_{n_1}}_{{}^{t}\wt{\bsym{f}}}=S^{{}^{t}\xi_{n_1}}_{{}^{t}\wt{\bsym{f}}^{[\bsym{\g}]}}$
 where ${}^{t}\xi_{n_1}:= {}^{t}\xi^{(1,\ldots,1)}_{n_1}$ and 
 ${}^{t}\wt{\bsym{f}}^{[\bsym{\g}]}=(\d^{0}f_{m+1}^{[\g_m]},\d^{1}f_{m}^{[\g_{m-1}]},\ldots,\d^{m}f_1^{[\g_0]})$.
 In fact, since the condition above means that
 $a_{\g_k}(n_1/e'_k)=a_{\g_k}(e'_{k+1}/e'_{k})$ for all $1\le k\le m-1$,
 the summand in \eqref{for:mid2} is written as 
\[
 \xi\bigl(\frac{n_1}{e'_1},\ldots,\frac{n_1}{e'_m}\bigr)f_1\bigl(e'_1\bigr)f_2^{[\g_1]}\bigl(\frac{e'_2}{e'_1}\bigr)\cdots
 f_m^{[\g_{m-1}]}\bigl(\frac{e'_m}{e'_{m-1}}\bigr)f_{m+1}^{[\g_m]}\bigl(\frac{n_1}{e'_m}\bigr)e'_1\cdots e'_m.
\]
 Hence, the claim follows.
\end{remark}

\begin{example}
 Retaining the notation in Example~\ref{ex:divisor}, one can easily see
 that
 ${}^{t}\wt{\bsym{f}}=(\d^{a_1+\cdots+a_m},\d^{a_1+\cdots+a_{m-1}+1},\ldots,\d^{a_1+m-1},\d^{m})$.
 Hence, the finite Fourier coefficient
 $\a_{\bsym{f},n_1}:=\a^{(1,\ldots,1),1_{m}}_{\bsym{f},n_1}$ of the 
 multiple divisor function $\s_{\bsym{a}}:=\s^{(1,\ldots,1)}_{\bsym{a}}$
 is given by  
\[
 \a_{\bsym{f},n_1}(d_2,\ldots,d_{m+1})
=n_1^{a_1+\cdots+a_m-m}\s_{\bsym{1}-{}^{t}\bsym{a}}\bigl(n_1,\frac{n_1}{d_{m+1}},\ldots,\frac{n_1}{d_{2}}\bigr),
\]
 where $\bsym{1}:=(1,\ldots,1)\in\bC^{m}$ and 
 ${}^{t}\bsym{a}:=(a_{m},\ldots,a_1)\in\bC^{m}$
\end{example}

\begin{example}
 Retaining the notation in Example~\ref{ex:gcd}, we have
 $S_{\bsym{f}}=f\circ \textrm{gcd}$. Again, by the degeneracy
 formula of \eqref{for:degeneracy2}, the finite Fourier
 coefficient $\a_{\bsym{f},n_1}$ of $f\circ \textrm{gcd}$ is given by
\begin{align*}
 \a_{\bsym{f},n_1}(d_2,\ldots,d_{m+1})
=n_1^{-m}S_{f*\m,\d^{m}}\bigl(n_1,\gcd{\frac{n_1}{d_2},\ldots,\frac{n_1}{d_{m+1}}}\bigr).
\end{align*}
\end{example}

\section{Dirichlet series attached to $S^{\bsym{\g}}_{\bsym{f}}$}
\label{sec:Dirichlet series}
 
 In this section, we examine both single-variable and multivariable
 Dirichlet series, the coefficients of which are given by the multiple
 Ramanujan sum $S^{\bsym{\g}}_{\bsym{f}}$.

\subsection{Single-variable Dirichlet series}

 We first examine a single-variable Dirichlet series. Recall the following 
 well-known formula concerning the Ramanujan sum $c(k,n)$ (see, e.g.,
 \cite{Titchmarsh1986}):  
\begin{equation}
\label{for:series}
 \sum^{\infty}_{k=1}c(k,n)k^{-s}=\frac{\s_{1-s}(n)}{\z(s)}
 \qquad (\Re(s)>1).
\end{equation}
 In this subsection, we give a generalization of this formula. For
 $j=1,2,\ldots,m+1$, let  
\[
 \Phi^{\bsym{\g}}_{\bsym{f}}(s;\check{\bsym{n}}_j)
:=\sum^{\infty}_{n_j=1}S^{\bsym{\g}}_{\bsym{f}}(n_1,\ldots,n_{m+1})n^{-s}_j,
\]
 where
 $\check{\bsym{n}}_j:=(n_1,\ldots,n_{j-1},n_{j+1},\ldots,n_{m+1})\in\bN^{m}$. 
 The following proposition says that the series
 $\Phi^{\bsym{\g}}_{\bsym{f}}(s;\check{\bsym{n}}_j)$ is written as
 the product of a Dirichlet series and a finite sum, which is again
 given by the multiple Ramanujan sum.

\begin{prop}
\label{prop:Dirichletfinite}

 $(\mathrm{i})$\ For $j=1$, we have  
\begin{equation}
\label{for:Phi1}
 \Phi^{\bsym{\g}}_{\bsym{f}}(s;\check{\bsym{n}}_1)
=L(s;f_1)S_{F_1}(n_2)=L(s;f_1)F_1(n_2) \qquad (\Re(s)>\s(f_1)),
\end{equation}
 where $F_1={F_{1,\bsym{f},\check{\bsym{n}}_1}^{\bsym{\g},s}}:=\d^{-s}\bigl(S^{(\g_{2},\ldots,\g_m)}_{f_2,\ldots,f_{m+1}}(\,\cdot\,,n_{3},\ldots,n_{m+1})*_{\g_{1}}\d^s\bigr)$.

 $(\mathrm{ii})$\ For $2\le j\le m+1$, we have
\begin{equation}
\label{for:Phij}
 \Phi^{\bsym{\g}}_{\bsym{f}}(s;\check{\bsym{n}}_j)
=\z(s)S^{(\g_1,\ldots,\g_{j-2})}_{f_1,\ldots,f_{j-2},F_{j}}(n_1,\ldots,n_{j-1}) \qquad (\Re(s)>1),
\end{equation}
 where $F_{j}={F_{j,\bsym{f},\check{\bsym{n}}_j}^{\bsym{\g},s}}:=\d^{-s}\bigl(S^{(\g_{j},\ldots,\g_m)}_{f_j,\ldots,f_{m+1}}(\,\cdot\,,n_{j+1},\ldots,n_{m+1})*_{\g_{j-1}}(\d^sf_{j-1})\bigr)$.
\end{prop}
\begin{proof}
 Suppose that $d^{\g_k}_k|\,\gcd{n_1,\ldots,n_{k+1}}$ for all
 $k=1,\ldots,m$. Then, $n_j$ is a multiple of
 $\lcm{d_{1}^{\g_{1}},\ldots,d_m^{\g_m}}=d_{1}^{\g_{1}}$ if $j=1$ and 
 $\lcm{d_{j-1}^{\g_{j-1}},\ldots,d_m^{\g_m}}=d_{j-1}^{\g_{j-1}}$ if
 $2\le j\le m+1$ (note that we only consider $m$-tuples
 $(d_1,\ldots,d_{m})$ such that $d_m^{\g_m}|\cdots |d_1^{\g_1}$).
 Hence, for $j=1$ and $\Re(s)>\s(f_1)$,
 $\Phi^{\bsym{\g}}_{\bsym{f}}(s;\check{\bsym{n}}_j)$ is expressed as
\[
 \sum_{d_k^{\g_k}|\,\gcd{n_2,\ldots,n_{k+1}}\atop
 (k=1,\ldots,m)}\sum^{\infty}_{l=1}f_{1}\bigl(\frac{d_1^{\g_1}l}{d_1^{\g_1}}\bigr)f_{2}\bigl(\frac{d_1^{\g_1}}{d_{2}^{\g_2}}\bigr)\cdots f_{m}\bigl(\frac{d_{m-1}^{\g_{m-1}}}{d_m^{\g_m}}\bigr)f_{m+1}\bigl(d_m^{\g_m}\bigr)\bigl(d_1^{\g_1}l\bigr)^{-s}\\
\]
 and, for $2\le j\le m+1$ and $\Re(s)>1$,
\begin{multline*}
 \sum_{d^{\g_k}_k|\,\gcd{n_1,\ldots,n_{k+1}}\atop
 (k=1,2,\ldots,j-2)}f_{1}\bigl(\frac{n_1}{d_1^{\g_1}}\bigr)
 f_{2}\bigl(\frac{d_{1}^{\g_{1}}}{d_{2}^{\g_{2}}}\bigr)\cdots
 f_{j-2}\bigl(\frac{d_{j-3}^{\g_{j-3}}}{d_{j-2}^{\g_{j-2}}}\bigr)\sum_{d_{j-1}^{\g_{j-1}}|d_{j-2}^{\g_{j-2}}}f_{j-1}\bigl(\frac{d_{j-2}^{\g_{j-2}}}{d_{j-1}^{\g_{j-1}}}\bigr)\\
 \times\sum_{d^{\g_k}_k|\,\gcd{d_{j-1}^{\g_{j-1}},n_{j+1},\ldots,n_{k+1}}\atop
 (k=j,\ldots,m)}\sum^{\infty}_{l=1}
 f_{j}\bigl(\frac{d_{j-1}^{\g_{j-1}}}{d_{j}^{\g_{j}}}\bigr)\cdots f_{m}\bigl(\frac{d_{m-1}^{\g_{m-1}}}{d_{m}^{\g_{m}}}\bigr)
 f_{m+1}\bigl(d_m^{\g_m}\bigr)\bigl(d_{j-1}^{\g_{j-1}}l\bigr)^{-s}.
\end{multline*}
 Thus, it is easy to see that these coincide, respectively, with
 \eqref{for:Phi1} and \eqref{for:Phij}. This completes the proof.
\end{proof}

\begin{example}
 For small $m$, the series
 $\Phi^{\bsym{\g}}_{\bsym{f}}(s;\check{\bsym{n}}_j)$ is explicitly given 
 as follows. For $m=1$, we have  
\begin{align}
\label{for:m=1j=1}
 \Phi^{(\g_1)}_{f_1,f_2}(s;n_2)
&=L(s;f_1)\sum_{d^{\g_1}|n_2}f_2(d^{\g_1})d^{-\g_1s},\\
\label{for:m=1j=2}
 \Phi^{(\g_1)}_{f_1,f_2}(s;n_1)
&=\z(s)\sum_{d^{\g_1}|n_1}f_1\bigl(\frac{n_1}{d^{\g_1}}\bigr)f_2(d^{\g_1})d^{-\g_1s}
\end{align}
 and for $m=2$, we have 
\begin{align*}
 \Phi^{(\g_1,\g_2)}_{f_1,f_2,f_3}(s;n_2,n_3)
&=L(s;f_1)\sum_{d_1^{\g_1}|n_2}\sum_{d_2^{\g_2}|\,\gcd{d_1^{\g_1},n_3}}f_2\bigl(\frac{d_1^{\g_1}}{d_2^{\g_2}}\bigr)f_3(d_2^{\g_2})d_1^{-\g_1s},\\
 \Phi^{(\g_1,\g_2)}_{f_1,f_2,f_3}(s;n_1,n_3)
&=\z(s)\sum_{d_1^{\g_1}|n_1}\sum_{d_2^{\g_2}|\,\gcd{d_1^{\g_1},n_3}}f_1\bigl(\frac{n_1}{d_1^{\g_1}}\bigr)f_2\bigl(\frac{d_1^{\g_1}}{d_2^{\g_2}}\bigr)f_3(d_2^{\g_2})d_1^{-\g_1s},\\
 \Phi^{(\g_1,\g_2)}_{f_1,f_2,f_3}(s;n_1,n_2)
&=\z(s)\sum_{d_1^{\g_1}|\,\gcd{n_1,n_2}}\sum_{d_2^{\g_2}|d_1^{\g_1}}f_1\bigl(\frac{n_1}{d_1^{\g_1}}\bigr)f_2\bigl(\frac{d_1^{\g_1}}{d_2^{\g_2}}\bigr)f_3(d_2^{\g_2})d_2^{-\g_2s}.
\end{align*}
 Formulas \eqref{for:m=1j=1} and \eqref{for:m=1j=2} were obtained in 
 \cite{AndersonApostol1953} for the case of $\g_1=1$. 
\end{example}

\begin{example}
\label{ex:geneRam}
 For $0\le k\le m+1$ and
 $\bsym{a}=(a_1,\ldots,a_{m+1-k})\in\bC^{m+1-k}$, we define
 generalizations of the classical Ramanujan sum $c(n_1,n_2)$ as
\[
 c^{\bsym{a}}_{m+1,k}(n_1,\ldots,n_{m+1})
:=S_{\underbrace{\m,\ldots,\m}_{k},\d^{a_1},\ldots,\d^{a_{m+1-k}}}(n_1,\ldots,n_{m+1}).
\] 
 It is clear that $c(n_1,n_2)=c^{(1)}_{2,1}(n_1,n_2)$. Then,
 for $1\le k\le m$, using formulas \eqref{for:Phi1} and
 $L(s;\mu)=1/\z(s)$, we can verify by induction on $m$ and $k$ that
\begin{multline}
\label{for:anbump}
 \sum^{\infty}_{n_1,\ldots,n_k=1}c^{\bsym{a}}_{m+1,k}(n_1,\ldots,n_{m+1})n_1^{-s_1}\cdots
 n_{k}^{-s_k}\\
=\frac{\prod^{k}_{j=2}\z(s_j)}{\prod^{k}_{j=1}\z(s_1+\cdots+s_j)}
\sum_{d|n_{k+1}}d^{a_1-(s_1+\cdots+s_{k})}\s_{\wt{\bsym{a}}}(d,n_{k+2},\ldots,n_{m+1}),
\end{multline}
 where
 $\wt{\bsym{a}}:=(a_2-a_1,a_3-a_2,\ldots,a_{m+1-k}-a_{m-k})\in\bC^{m-k}$
 (see Example~\ref{ex:divisor}). Setting $m=1$, $k=1$ (in this case,
 $\s_{\wt{\bsym{a}}}\equiv 1$) and $a_1=1$, we obtain formula
 \eqref{for:series}. See \cite{Bump1984} for a multivariable analogue of
 formula \eqref{for:series}.
\end{example}

\begin{example}
 Retaining the notation in Example~\ref{ex:gcd}, we have
 $S_{\bsym{f}}=f\circ \textrm{gcd}$. Let 
 $\gcd{\check{\bsym{n}}_j}:=\gcd{n_1,\ldots,n_{j-1},n_{j+1},\ldots,n_{m+1}}$. 
 Then, since it holds that
 $\gcd{n_1,\ldots,n_{m+1}}=\gcd{n_j,\gcd{\check{\bsym{n}}_j}}$,
 we have the following well-known formula
\begin{align*}
 \sum^{\infty}_{n_j=1}(f\circ \textrm{gcd})(n_1,\ldots,n_{m+1})n_j^{-s}
&=\Phi_{\d^{0},f*\m}(s;\gcd{\check{\bsym{n}}_j})=\z(s)\sum_{d|\,\gcd{\check{\bsym{n}}_j}}(f*\m)(d)d^{-s}.
\end{align*}  
\end{example}

\subsection{Multivariable Dirichlet series}

 For an arithmetic function $F(n_1,\ldots,n_k)$ of $k$-variables, we
 denote $L(\bsym{s};F)$ by the multivariable Dirichlet series attached
 to $F$, that is,  
\[
 L(\bsym{s};F)
:=\sum^{\infty}_{n_1,\ldots,n_k=1}F(n_1,\ldots,n_k)n_1^{-s_1}\cdots
n_{k}^{-s_k}
\]
 for $\bsym{s}=(s_1,\ldots,s_k)\in\Sigma(F)$, where
 $\Sigma(F)\subset \bC^{k}$ is the region of absolute convergence for  
 $L(\bsym{s};F)$. The first goal of this subsection is to calculate
 $L(\bsym{s};S^{\bsym{\g}}_{\bsym{f}})$ explicitly.

\begin{thm}
\label{thm:DirichletS}
 Let $\g_0|\g_1|\g_2|\cdots|\g_m$ with $\g_0\in\bN$. Define the function 
 $a_{\g_0}S^{\bsym{\g}}_{\bsym{f}}$ of $(m+1)$-variables as 
 $(a_{\g_0}S^{\bsym{\g}}_{\bsym{f}})(n_1,\ldots,n_{m+1})$
 $:=a_{\g_0}(n_1)S^{\bsym{\g}}_{\bsym{f}}(n_1,\ldots,n_{m+1})$. Then, we
 have   
\begin{equation}
\label{for:DirichletSS}
 L(\bsym{s};a_{\g_0}S^{\bsym{\g}}_{\bsym{f}})
=\prod^{m+1}_{j=2}\z(s_j)\cdot
 \prod^{m+1}_{j=1}L(s_1+\cdots+s_j;f_j^{[\g_{j-1}]}) 
\end{equation}
 in the region 
\[
\left\{
 \bsym{s}=(s_1,\ldots,s_{m+1})\in\bC^{m+1}
\left|
\begin{array}{ll}
 \Re(s_j)>1 & (2\le j\le m+1),\\[3pt]
 \Re(s_1+\cdots+s_j)>\s(f_j) & (1\le j\le m+1)
\end{array}
\right.
\right\}.
\]
 In particular, setting $\g_0=1$, we have
\begin{equation}
\label{for:DirichletS} 
 L(\bsym{s};S^{\bsym{\g}}_{\bsym{f}})
=\prod^{m+1}_{j=2}\z(s_j)\cdot L(s_1;f_1)\prod^{m+1}_{j=2}L\bigl(s_1+\cdots+s_j;f_j^{[\g_{j-1}]}\bigr).
\end{equation}
\end{thm}
\begin{proof}
 This is proven by induction on $m$. Suppose that $m=1$. Then, it holds
 from formula \eqref{for:m=1j=2} that
\begin{align*}
 L\bigl((s_1,s_2);a_{\g_0}S^{(\g_1)}_{f_1,f_2}\bigr)
&=\sum^{\infty}_{n_1=1}\Phi^{(\g_1)}_{f_1,f_2}(s_2;n_1^{\g_0})n_1^{-\g_0 s}\\
&=\z(s_2)\sum^{\infty}_{n_1=1}\sum_{d^{\g_1}|n_1^{\g_0}}f_1\bigl(\frac{n_1}{d^{\g_1}}\bigr)f_2\bigl(d^{\g_1}\bigr)d^{-\g_1 s_2}n_1^{-\g_0 s_1}.
\end{align*}
 Here, since $\g_0|\g_1$, $n_1^{\g_0}$ is expressed as
 $n_1^{\g_0}=d^{\g_1}l^{\g_0}$ with $l\in\bN$. Hence, this expression is
 written as
\begin{align*}
 \z(s_2)\sum^{\infty}_{d=1}\sum^{\infty}_{l=1}f_1\bigl(\frac{d^{\g_1}l^{\g_0}}{d^{\g_1}}\bigr)f_2\bigl(d^{\g_1}\bigr)d^{-\g_1
 s_2}\bigl(d^{\g_1}l^{\g_0}\bigr)^{-s_1},
\end{align*}
 thus completing the proof for $m=1$. Next, suppose that $m-1$. Note
 that   
 $S^{\bsym{\g}}_{\bsym{f}}(n_1,\ldots,n_{m+1})$ $=S^{(\g_1)}_{f_1,S^{(\g_2,\ldots,\g_m)}_{f_2,\ldots,f_{m+1}}(\,\cdot\,,n_3,\ldots,n_{m+1})}(n_1,n_2)$
 from formula \eqref{for:Srec}. Then, from \eqref{for:DirichletSS} for
 $m=1$, $L(\bsym{s};a_{\g_0}S^{\bsym{\g}}_{\bsym{f}})$ is equal to  
\begin{align*}
&\sum^{\infty}_{n_3,\ldots,n_{m+1}=1}L\bigl((s_1,s_2);a_{\g_0}S^{(\g_1)}_{f_1,S^{(\g_2,\ldots,\g_m)}_{f_2,\ldots,f_{m+1}}(\,\cdot\,,n_3,\ldots,n_{m+1})}\bigr)n_3^{-s_1}\cdots
 n_{m+1}^{-s_{m+1}}\\
&=\z(s_2)L(s_1;f^{[\g_0]}_1)\sum^{\infty}_{n_3,\ldots,n_{m+1}=1}L\bigl(s_1+s_2;a_{\g_1}S^{(\g_2,\ldots,\g_m)}_{f_2,\ldots,f_{m+1}}(\,\cdot\,,n_3,\ldots,n_{m+1})\bigr)n_3^{-s_1}\cdots
 n_{m+1}^{-s_{m+1}}\\
&=\z(s_2)L(s_1;f^{[\g_0]}_1)L\bigl((s_1+s_2,s_3,\ldots,s_{m+1});a_{\g_1}S^{(\g_2,\ldots,\g_m)}_{f_2,\ldots,f_{m+1}}\bigr).
\end{align*}
 This completes the proof for $m$ based on the assumption of
 induction. Hence, we obtain the desired formula.
\end{proof}
 
 Remark that formula \eqref{for:Dirichletseries} is regarded as the 
 ``diagonally-summed version'' of \eqref{for:DirichletS} since
 $S^{(\g_1,\ldots,\g_m)}_{f_1,\ldots,f_{m+1}}(n,\ldots,n)=(f_{m+1}*_{\g_m}\cdots*_{\g_1}f_1)(n)$. 

\begin{example}
 Retaining the notation in Example~\ref{ex:divisor}, we have
\begin{equation}
\label{for:divisorL}
 L(\bsym{s};\s_{\bsym{a}})
=\z(s_1)\prod^{m+1}_{k=2}\z(s_j)\z(s_1+\cdots+s_{j}-a_1-\cdots-a_{j-1}).
\end{equation}
\end{example}

\begin{example}
 Retaining the notation in Example~\ref{ex:geneRam}, we have 
\[
 L\bigl(\bsym{s};c^{\bsym{a}}_{m+1,k}\bigr)
=\frac{\prod^{m+1}_{l=2}\z(s_{l})\prod^{m+1-k}_{l=1}\z(s_1+\cdots+s_{k+l}-a_l)}{\prod^{k}_{l=1}\z(s_1+\cdots+s_l)}.
\]
 We can also obtain this formula from \eqref{for:anbump} and
 \eqref{for:divisorL}. In particular, we have 
 $L((s_1,s_2);c)=\z(s_2)\z(s_1+s_2-1)/\z(s_1)$.
\end{example}

\begin{example}
 Retaining the notation in Example~\ref{ex:gcd}, we have
 $S_{\bsym{f}}=f\circ\textrm{gcd}$. In this case, it holds that 
\[
 L\bigl(\bsym{s};f\circ\textrm{gcd}\bigr)
=\frac{\z(s_1)\cdots \z(s_{m+1})}{\z(s_1+\cdots+s_{m+1})}
L(s_1+\cdots+s_{m+1};f).
\] 
\end{example}

 Next, we present another type of multivariable Dirichlet series
 involving $S^{\bsym{\g}}_{\bsym{f}}$. Recall the well-known Ramanujan
 formula for the divisor function $\s_a$ (see e.g.,
 \cite{Titchmarsh1986}, also \cite{Bump1992}):
\[
 L(s;\s_a\s_b)
=\frac{\z(s)\z(s-a)\z(s-b)\z(s-a-b)}{\z(2s-a-b)}
\] 
 for $\Re(s)>\max\{1,\Re(a)+1,\Re(b)+1,\Re(a+b)+1\}$.
 This formula was generalized by J. M. Borwein and K. K. Choi
 \cite{BorweinChoi2003} as follows. Let $f_1,f_2,g_1,g_2\in\cM^{c}$. 
 Then, we have    
\begin{equation}
\label{for:didouble}
 L\bigl(s;(f_2*_{\g}f_1)(g_2*_{\g}g_1)\bigr)
=\frac{L\bigl(s;f_1g_1\bigr)L\bigl(s;(f_2g_1)^{[\g]}\bigr)L\bigl(s;(f_1g_2)^{[\g]}\bigr)L\bigl(s;(f_2g_2)^{[\g]}\bigr)}{L\bigl(2s;(f_1f_2g_1g_2)^{[\g]}\bigr)}\quad
\end{equation} 
 for
 $\Re(s)>\max\{\s(f_1g_1),\s(f_2g_1),\s(f_1g_2),\s(f_2g_2),\s(f_1f_2g_1g_2)/2\}$.
 They gave formula \eqref{for:didouble} for the case in which $\g=1$,
 and it is easy to obtain the equation for general $\g\in\bN$ in the
 same manner. Similarly, since
 $S^{(\g)}_{f_1,f_2}(n,n)=(f_2*_{\g}f_1)(n)$, we can regard formula
 \eqref{for:didouble} as the diagonally-summed version of the following
 formula.

\begin{thm}
 Let $f_1,f_2,g_1,g_2\in\cM^{c}$. Then, we have
\begin{multline}
\label{for:double}
 L\bigl((s_1,s_2);S^{(\g)}_{f_1,f_2} S^{(\g)}_{g_1,g_2}\bigr)\\
=\frac{\z(s_2)L\bigl(s_1;f_1g_1\bigr)L\bigl(s_1+s_2;(f_2g_1)^{[\g]}\bigr)L\bigl(s_1+s_2;(f_1g_2)^{[\g]}\bigr)L\bigl(s_1+s_2;(f_2g_2)^{[\g]}\bigr)}{L\bigl(2(s_1+s_2);(f_1f_2g_1g_2)^{[\g]}\bigr)}
\end{multline}
 in the region
\[
\left\{(s_1,s_2)\in\bC^2
\left|
\begin{array}{l}
 \Re(s_2)>1,\ \Re(s_1)>\s(f_1g_1),\\[3pt]
 \Re(s_1+s_2)>\max\{\s(f_2g_1),\s(f_1g_2),\s(f_2g_2),\s(f_1f_2g_1g_2)/2\}
\end{array}
\right.
\right\}.
\]
\end{thm} 
\begin{proof}
 First, note that, for a multiplicative arithmetic function $F$ of
 $k$-variables, $L(\bsym{s};F)$ has the Euler product expression
 $L(\bsym{s};F)=\prod_{p}\sum^{\infty}_{l_1,\ldots,l_{k=1}}F(p^{l_1},\ldots,p^{l_k})p^{-s_1l_1}\cdots p^{-s_kl_k}$.

 Let $h_1,h_2\in\cM^{c}$. Then, it holds that  
\[
 S^{(\g)}_{h_1,h_2}(p^{l_1},p^{l_2})=\frac{h_1(p)^{\g}}{h_1(p)^{\g}-h_2(p)^{\g}}h_1(p)^{l_1}\Bigl(1-\bigl(\frac{h_2(p)}{h_1(p)}\bigr)^{\g\cdot m^{\g}(l_1,l_2)}\Bigr),
\]
 where $m^{\g}(l_1,l_2):=\Gauss{\frac{1}{\g}\min\{l_1,l_2\}}+1$.
 Hence, from the Euler product expression, we have
\begin{align}
\label{for:MMM}
& L\bigl((s_1,s_2);S^{(\g)}_{f_1,f_2} S^{(\g)}_{g_1,g_2}\bigr)\\
&=\prod_p\sum^{\infty}_{l_1,l_2=0}S^{(\g)}_{f_1,f_2}(p^{l_1},p^{l_2})S^{(\g)}_{g_1,g_2}(p^{l_1},p^{l_2})p^{-s_1l_1}p^{-s_2l_s}\nonumber\\
&=\prod_p \bigl(C^{\g}_p\bigr)^{-1}
\sum^{\infty}_{l_1,l_2=0}(a_1b_1x_1)^{l_1}x_2^{l_2}\Bigl(1-\bigl(\frac{a_2}{a_1}\bigr)^{\g\cdot m^{\g}(l_1,l_2)}\Bigr)\Bigl(1-\bigl(\frac{b_2}{b_1}\bigr)^{\g\cdot m^{\g}(l_1,l_2)}\Bigr),\nonumber
\end{align}
 where $a_i:=f_i(p)$, $b_i:=g_i(p)$, $x_i:=p^{-s_i}$ for $i=1,2$ and
 $C^{\g}_p=C^{\g}_p(f_1,f_2,g_1,g_2):=(a_1^{\g}-a_2^{\g})(b_1^{\g}-b_2^{\g})/(a_1b_1)^{\g}$
 for each prime $p$. Write $I$ for the inner sum of the right-most hand
 side of \eqref{for:MMM}. Moreover, divide $I$ into two parts as
 $\sum^{\infty}_{l_1,l_2=0}=\sum_{l_1\ge l_2}+\sum_{l_1<l_2}=\sum^{\infty}_{l_2=0}\sum^{\infty}_{l_1=l_2}+\sum^{\infty}_{l_1=0}\sum^{\infty}_{l_2=l_1+1}$
 and denote $I_1$ and $I_2$ by the former and latter sums,
 respectively. Writing $l_2$ as $l_2=\g l+k$ in $I_1$ and $l_1$ as
 $l_1=\g l+k$ in $I_2$ with $l\in\bZ_{\ge 0}$ and $0\le k\le \g-1$,
 respectively, we have 
\begin{align*}
 I_1
&=\sum^{\g-1}_{k=0}\sum^{\infty}_{l=0}\sum^{\infty}_{l_1=\g
 l+k}(a_1b_1x_1)^{l_1}x_2^{\g l+k}\Bigl(1-\bigl(\frac{a_2}{a_1}\bigr)^{\g(l+1)}\Bigr)\Bigl(1-\bigl(\frac{b_2}{b_1}\bigr)^{\g(l+1)}\Bigr)\\
&=\frac{1}{1-a_1b_1x_1}\frac{1-(a_1b_1x_1x_2)^{\g}}{1-a_1b_1x_1x_2}M
\intertext{and similarly}
I_2
&=\frac{x_2}{1-x_2}\frac{1-(a_1b_1x_1x_2)^{\g}}{1-a_1b_1x_1x_2}M,
\end{align*}
 where
 $M:=\sum^{\infty}_{l=0}(a_1b_1x_1x_2)^{\g l}\bigl(1-(\frac{a_2}{a_1})^{\g(l+1)}\bigr)\bigl(1-(\frac{b_2}{b_1})^{\g(l+1)}\bigr)$.
 Hence, we have
\[
 I
=I_1+I_2
=\frac{1-(a_1b_1x_1x_2)^{\g}}{(1-x_2)(1-a_1b_1x_1)}M.
\]
 Now, $M$ is straightforwardly calculated as 
\[
 \frac{C^{\g}_p\bigl(1-(a_1a_2b_1b_2x_1^2x_2^2)^{\g}\bigr)}{\bigl(1-(a_1b_1x_1x_2)^{\g}\bigr)\bigl(1-(a_2b_1x_1x_2)^{\g}\bigr)\bigl(1-(a_1b_2x_1x_2)^{\g}\bigr)\bigl(1-(a_2b_2x_1x_2)^{\g}\bigr)},
\]
 where the desired formula immediately follows from formula
 \eqref{for:MMM}.
\end{proof}

\begin{remark}
 For $m\ge 2$, we cannot expect that the single-variable Dirichlet
 series
 $L(s;(f_{m+1}*_{\g_m}\cdots*_{\g_1}f_1)(g_{m+1}*_{\g_m}\cdots*_{\g_1}g_1))$
 has a product expression such as \eqref{for:didouble} and, similarly,
 the multivariable Dirichlet series
 $L\bigl((s_1,\ldots,s_{m+1});S^{(\g_1,\ldots,\g_m)}_{f_1,\ldots,f_{m+1}}S^{(\g_1,\ldots,\g_m)}_{g_1,\ldots,g_{m+1}}\bigr)$
 has a product expression such as \eqref{for:double}.
\end{remark}

\section{A generalization of the Smith determinant}
\label{sec:determinant}

 In this section, we evaluate the hyperdeterminants in the sense of
 Cayley of the hypermatrices, the entries of which are given by values
 of the multiple Ramanujan sum $S^{\bsym{\g},\xi}_{\bsym{f}}$. 

\subsection{Hyperdeterminant}

 Recall the definition and some properties of the hyperdeterminants in
 the sense of Cayley. Let
 $A=(A(i_1,\ldots,i_k))_{1\le i_1,\ldots,i_k\le n}$ be a $k$-dimensional
 matrix of order $n$. For a subset $I\subseteq\{1,\ldots,k\}$, set
 $\e_j=1$ if $j\in I$, and $0$ otherwise. Then, the hyperdeterminant
 $\det_{I}A$ of $A$ with the signature $I$ is defined as
\begin{equation}
\label{def:hyperdeterminant}
 \det_{I}A
:=\frac{1}{n!}\sum_{\s_1,\ldots,\s_k\in\fS_{n}}\prod^{k}_{j=1}\sgn(\s_j)^{\e_j}\sum^{n}_{v=1}A(\s_1(v),\ldots,\s_k(v)),
\end{equation}
 where $\fS_n$ is the symmetric group of degree $n$. We obtain the usual
 determinant of a matrix $A$ when $k=2$ and $I=\{1,2\}$;
 $\det_{\{1,2\}}A=\det\,{A}$. Note that $\det_{I}A$ is identically zero
 if the number of elements in $I$ is odd. For the theory of 
 hyperdeterminants, see, e.g., \cite{Oldenburger1940}. We next present
 some lemmas that will be needed in the subsequent discussion. 

\begin{lem}
 Let $A=(A(i_1,\ldots,i_k))$ be a $k$-dimensional matrix of order $n$. 

 $(\mathrm{i})$\ For $\pi\in\fS_n$, we have 
\begin{equation}
\label{for:s}
 \det_{I}\bigl(A(\pi(i_1),\ldots,\pi(i_k))\bigl)
=\det_{I}\bigl(A(i_{1},\ldots,i_{k})\bigl).
\end{equation}

 $(\mathrm{ii})$\ For $\pi\in\fS_k$, set
 ${\pi}^{-1}(I):=\{{\pi}^{-1}(i)\,|\,i\in I\}$. Then, we have  
\begin{equation}
\label{for:t}
 \det_{I}\bigl(A(i_{\pi(1)},\ldots,i_{\pi(k)})\bigl)
=\det_{{\pi}^{-1}(I)}\bigl(A(i_{1},\ldots,i_{k})\bigl).
\end{equation}
\end{lem}
\begin{proof}
 These are immediately obtained from the definition of the
 hyperdeterminant. 
\end{proof}

 Let $A=(A(i_1,\ldots,i_{k}))$ and $B=(B(i_1,\ldots,i_l))$ be
 $k$-dimensional and $l$-dimensional matrices, respectively, of order
 $n$. Then, the Cayley product $AB$ of $A$ and $B$ is a
 $(k+l-2)$-dimensional matrix of order $n$ given by 
\begin{equation}
\label{def:Cayley product}
 (AB)(i_1,\ldots,i_{k+l-2})
:=\sum^{n}_{j=1}A(i_1,\ldots,i_{k-1},j)B(j,i_{k},\ldots,i_{k+l-2}).
\end{equation}

\begin{lem}
 Let $A=(A(i_1,\ldots,i_k))$ and $B=(B(i_1,\ldots,i_l))$ be 
 $k$-dimensional and $l$-dimensional matrices, respectively, of order
 $n$. Let $K\subseteq \{1,2,\ldots,k-1\}$ and
 $L\subseteq\{2,3,\ldots,l\}$ with odd cardinality. Set
 $I:=K\cup (L+(k-2))$, where $L+(k-2):=\{l+k-2\,|\,l\in L\}$. Then, we
 have 
\begin{equation}
\label{for:product}
 \det_{I}(AB)=\det_{K\cup\{k\}}A\cdot \det_{\{1\}\cup L}B.
\end{equation}  
\end{lem}
\begin{proof}
 See \cite{Haukkanen1992,Rice1930}.
\end{proof}

\subsection{Smith hyperdeterminants}

 Let $S=\{x_1,\ldots, x_n\}$ be a set of distinct positive integers. The 
 goal of this subsection is to evaluate the hyperdeterminant 
 $\det_{I}(S^{\bsym{\g},\xi}_{\bsym{f}}(x_{i_1},\ldots,x_{i_{m+1}}))$
 for a factor-closed set $S$. Here, $S$ is called factor-closed if $S$
 contains every divisor of $x$ for any $x\in S$. To accomplish this, we
 slightly extend the definition of the even arithmetic function examined
 in Section~\ref{sec:definition}.

 First, we define some notations. For $\bsym{d}=(d_1,\ldots,d_k)\in\bN^{k}$,
 we set $|\bsym{d}|:=\sum^{k}_{j=1}d_j$ and
 $x_{\bsym{d}}:=(x_{d_1},\ldots,x_{d_k})$. Moreover, for $r\in\bN$, we
 set $\gcd{\bsym{d},r}:=(\gcd{d_{1},r},\ldots,\gcd{d_{k},r})\in\bN^{k}$,
 write $\bsym{d}|r$ if $d_j|r$ for $1\le j\le k$ and, in this case,
 set $r/\bsym{d}:=(r/d_{1},\ldots,r/d_{k})\in\bN^{k}$. Let $f(m,n)$ be
 an arithmetic function of $2$-variables. For
 $\bsym{m}=(m_1,\ldots,m_k)\in\bN^{k}$ and
 $\bsym{n}=(n_1,\ldots,n_k)\in\bN^{k}$, we set
 $f(\bsym{n},\bsym{m}):=\prod^{k}_{j=1}f(m_j,n_j)$. 

 Let $\bsym{r}=(r_1,\ldots,r_m)\in\bN^{m}$ and 
 $\bsym{k}=(k_1,\ldots,k_m)\in\bN^{m}$. Let
 $\bsym{n}_j=(n_{j,1},\ldots,n_{j,k_j})\in\bN^{k_j}$ for
 $j=1,\ldots,m$. We call an arithmetic function
 $F(\bsym{r};\bsym{n}_1,\ldots,\bsym{n}_m)$ of $|\bsym{k}|$-variables
 $\bsym{n}_1,\ldots,\bsym{n}_m$ even $\pmod{\bsym{r}^{(\bsym{k})}}$ if
 it holds that  
\[
 F\bigl(\bsym{r};\gcd{\bsym{n}_1,r_1},\ldots,\gcd{\bsym{n}_m,r_m}\bigr)
=F(\bsym{r};\bsym{n}_1,\ldots,\bsym{n}_m).
\]
 Then, the following lemma is obtained as described in \cite{Cohen1960}
 (note that the case of $m=1$ is nothing more than Theorem~$22$ in
 \cite{Cohen1960}).

\begin{lem}
 A function $F(\bsym{r};\bsym{n}_1,\ldots,\bsym{n}_m)$ is even
 $\pmod{\bsym{r}^{(\bsym{k})}}$ if and only if it possesses a
 representation of the form
\begin{equation}
\label{for:generaleven}
 F(\bsym{r};\bsym{n}_1,\ldots,\bsym{n}_m)
=\sum_{\bsym{d}_j|r_j \atop (j=1,\ldots,m)}\a_{\bsym{r}}(\bsym{d}_1,\ldots,\bsym{d}_{m})\prod^{m}_{j=1}c(\bsym{d}_j,\bsym{n}_j),
\end{equation}
 where the sum is taken over all
 $\bsym{d}_j=(d_{j,1},\ldots,d_{j,k_j})\in\bN^{k_j}$ for $j=1,\ldots,m$
 satisfying $\bsym{d}_{j}|r_j$. Moreover, the finite Fourier
 coefficients $\a_{\bsym{r}}(\bsym{d}_1,\ldots,\bsym{d}_{m})$ 
 are uniquely determined by 
\begin{equation}
\label{for:generalevencoeff}
 \a_{\bsym{r}}(\bsym{d}_1,\ldots,\bsym{d}_{m})
=\frac{1}{r^{k_1}_1\cdots r^{k_m}_m}
\sum_{\bsym{\d}_j|r_j \atop (j=1,\ldots,m)}
F\bigl(\bsym{r};\bsym{\d}_1,\ldots,\bsym{\d}_m\bigr)
\prod^{m}_{j=1}
c\bigl(\frac{r_j}{\bsym{\d}_j},\frac{r_j}{\bsym{d}_j}\bigr),
\end{equation}
 where the sum is also taken over all 
 $\bsym{\d}_j=(\d_{j,1},\ldots,\d_{j,k_j})\in\bN^{k_j}$ for $j=1,\ldots,m$ 
 satisfying $\bsym{\d}_{j}|r_j$. \qed 
\end{lem} 

 Now, we obtain the following proposition.

\begin{prop}
\label{prop:hypdet}
 Let $S=\{x_1,\ldots,x_n\}$ be a factor-closed set, and let
 $F(\bsym{r};\bsym{n}_1,\ldots,\bsym{n}_k)$ be an even function
 $\pmod{\bsym{r}^{(\bsym{k})}}$. Define the
 $(m+|\bsym{k}|)$-dimensional matrix $B$ of order $n$ as 
\[
 B(i_1,\ldots,i_{m+|\bsym{k}|})
=F\bigl(x_{i_1},\ldots,x_{i_m};\underbrace{x_{i_{m+1}},\ldots,x_{i_{m+k_1}}}_{k_1},\ldots,\underbrace{x_{i_{m+k_1+\ldots+k_{m-1}+1}},\ldots,x_{i_{m+|\bsym{k}|}}}_{k_m}\bigr).
\]
 Let $I$ be a subset of $\{1,2,\ldots,m+|\bsym{k}|\}$ with even
 cardinality such that $\{m+1,m+2\ldots,m+|\bsym{k}|\}\subseteq I$, and
 set $\e_j=1$ if $j\in I$, and $0$ otherwise. Define the subset $\wt{I}$
 of $\{1,2,\ldots,m\}$ as $j\in \wt{I}$ if and only if $\e_j+k_j$ is odd
 for $1\le j\le m$. Then, we have
\begin{equation}
\label{for:hypdetm}
 \det_{I}B
=(x_1\cdots
x_n)^{|\bsym{k}|}\det_{\widetilde{I}}\bigl(\a_{x_{i_1},\ldots,x_{i_{m}}}(\underbrace{x_{i_1},\ldots,x_{i_1}}_{k_1},\ldots,\underbrace{x_{i_m},\ldots,x_{i_m}}_{k_m})\bigr)_{1\le
i_1,\ldots,i_{m}\le n}.
\end{equation}
 Here, $\a_{x_{i_1},\ldots,x_{i_{m}}}$ is the finite Fourier coefficient
 of $F$ given by \eqref{for:generalevencoeff}.
\end{prop}

 We then need the following lemma. 

\begin{lem}
 Let $A=(A(i_1,\ldots,i_k))$ and $C=(C(i_1,i_2))$ be $k$-dimensional and
 $2$-dimensional matrices, respectively, of order $n$. For
 $l=0,1,\ldots,k$, define the $k$-dimensional matrices $A^{(l)}_C$ of
 order $n$ by the following recursion formula:
\[
\left\{
\begin{array}{ll}
 A^{(0)}_{C}(i_1,\ldots,i_k):=A(i_1,\ldots,i_k) &\ \ (l=0),\\[3pt]
 \,A^{(l)}_{C}(i_1,\ldots,i_k):=(A^{(l-1)}_{C}C)(i_2,\ldots,i_k,i_1)
 &\ \ (l=1,2,\ldots,k),
\end{array}
\right.  
\]
 where $A^{(l-1)}_CC$ is the Cayley product of $A^{(l-1)}_C$ and
 $C$. Then, we have 
\begin{equation}
\label{for:recA}
 A^{(l)}_{C}(i_1,\ldots,i_k)
=\sum^{n}_{j_1,\ldots,j_l=1}
A(i_{l+1},\ldots,i_k,j_1,\ldots,j_{l})\prod^{l}_{h=1}C(j_h,i_h).
\end{equation}
\end{lem}
\begin{proof}
 This is shown by induction on $l$ from the definition of the Cayley
 product \eqref{def:Cayley product}. 
\end{proof}

\begin{proof}
[Proof of Proposition~\ref{prop:hypdet}]
 From formula \eqref{for:s}, we can assume that
 $x_1<x_2<\cdots<x_n$. Define the $(m+|\bsym{k}|)$-dimensional matrix
 $A$ and $2$-dimensional matrices $C$ of order $n$, respectively, as 
\begin{align*}
 A(i_1,\ldots,i_{m+|\bsym{k}|})
:&=\a_{x_{i_1},\ldots,x_{i_{m}}}(x_{i_{m+1}},\ldots,x_{i_{m+|\bsym{k}|}})\prod^{m}_{j=1}\prod^{k_j}_{l_j=1}e_{i_{m+k_1+\cdots+k_{j-1}+l_j},i_j},\\
 C(i_1,i_2)
:&=c(x_{i_1},x_{i_2}),
\end{align*}
 where $e_{j,i}=1$ if $x_j|x_i$, and $0$ otherwise. Let
 $\r=(1,2,\cdots,m+|\bsym{k}|)\in\fS_{m+|\bsym{k}|}$ be the cyclic
 permutation of order $(m+|\bsym{k}|)$. Then, we have
\begin{equation}
\label{for:B}
 B(i_{1},\ldots,i_{m+|\bsym{k}|})
=A^{(|\bsym{k}|)}_{C}(i_{\r^{m}(1)},\ldots,i_{\r^{m}(m+|\bsym{k}|)}).
\end{equation}
 In fact, by equation \eqref{for:recA}, the right-hand side of
 \eqref{for:B} is equal to  
\begin{align*}
 A^{(|\bsym{k}|)}_{C}&(i_{m+1},\ldots,i_{m+|\bsym{k}|},i_1,\ldots,i_{m})\\
&=\sum^{n}_{d_{j,1},\ldots,d_{j,k_j}=1 \atop
 (j=1,\ldots,m)}\a_{x_{i_1},\ldots,x_{i_{m}}}\bigl(\underbrace{x_{d_{1,1}},\ldots,x_{d_{1,k_1}}}_{k_1},\ldots,\underbrace{x_{d_{m,1}},\ldots,x_{d_{m,k_m}}}_{k_m}\bigr)\\
&\ \ \ \times\prod^{m}_{j=1}\prod^{k_j}_{l_j=1}e_{d_{j,l_j},i_j}c\bigl(x_{d_{j,l_j}},x_{i_{m+k_1+\cdots+k_{j-1}+l_j}}\bigr)\\
&=\sum_{x_{\bsym{d}_j}|x_{i_j} \atop (j=1,\ldots,m)}\a_{x_{i_1},\ldots,x_{i_{m}}}(x_{\bsym{d}_1},\ldots,x_{\bsym{d}_{m}})\prod^{m}_{j=1}\prod^{k_j}_{l_j=1}c\bigl(x_{d_{j,l_j}},x_{i_{m+k_1+\cdots+k_{j-1}+l_j}}\bigr),
\end{align*}
 where $\bsym{d}_j:=(d_{j,1},\ldots,d_{j,k_j})$ for $1\le j\le m$.
 This clearly coincides with the left-hand side of formula \eqref{for:B}
 from the expression of \eqref{for:generaleven} of $F$, thus completing
 the proof. Note that, since $S$ is factor-closed, every divisor of
 $x_j$ can again be written as $x_d$ for some $d$. Next, using formula
 \eqref{for:t}, we have
\begin{align*}
 \det_{I}B
&=\det_{I}\bigl(A^{(|\bsym{k}|)}_{C}(i_{\r^{m}(1)},\ldots,i_{\r^{m}(m+|\bsym{k}|)})\bigr)\\
&=\det_{I}\bigl((A^{(|\bsym{k}|-1)}_{C}C)(i_{\r^{m}(2)},\ldots,i_{\r^{m}(m+|\bsym{k}|)},i_{\r^{m}(1)})\bigr)\\
&=\det_{I}\bigl((A^{(|\bsym{k}|-1)}_{C}C)(i_{\r^{m+1}(1)},\ldots,i_{\r^{m+1}(m+|\bsym{k}|)})\bigr)\\
&=\det_{\r^{-m-1}(I)}\bigl((A^{(|\bsym{k}|-1)}_{C}C)(i_1,\ldots,i_{m+|\bsym{k}|})\bigr).
\end{align*}
 Here, note that $m+1\in I$ implies
 $m+|\bsym{k}|=\r^{-m-1}(m+1)\in \r^{-m-1}(I)$. Therefore, using formula
 \eqref{for:product} with $K=\r^{-m-1}(I)\setminus \{m+|\bsym{k}|\}$ and
 $L=\{2\}$, we see that the right-hand side of the above equation is
 equivalent to    
\begin{align*}
&\det_{(\r^{-m-1}(I)\setminus
 \{m+|\bsym{k}|\})\cup\{m+|\bsym{k}|\}}\bigl(A^{(|\bsym{k}|-1)}_{C}\bigr)\cdot
 \det_{\{1,2\}}C\\
&\ \ =\det\,{C} \cdot
 \det_{\r^{-m-1}(I)}\bigl((A^{(|\bsym{k}|-2)}_{C}C)(i_2,\ldots,i_{m+|\bsym{k}|},i_1)\bigr)\\
&\ \ =\det\,{C} \cdot
 \det_{\r^{-m-2}(I)}\bigl((A^{(|\bsym{k}|-2)}_{C}C)(i_1,\ldots,i_{m+|\bsym{k}|})\bigr).
\end{align*}
 Under the condition $m+1,\ldots,m+|\bsym{k}|\in I$, the above procedure
 yields  
\begin{align*}
 \det_{I}B
&=\bigl(\det_{\{1,2\}}C\bigr)^{|\bsym{k}|} \cdot
 \det_{\r^{-m-|\bsym{k}|}(I)}\bigl(A^{(0)}_{C}\bigr)=\bigl(\det_{\{1,2\}}C\bigr)^{|\bsym{k}|} \cdot
 \det_{I}A.
\end{align*}
 By Corollary~$3$ in \cite{BourqueLigh1993}, we have
 $\det\,{C}=\det\bigl(c(x_{i_1},x_{i_2})\bigr)=x_1\cdots x_n$.
 Hence, it is sufficient to calculate $\det_{I}A$. Since $\e_j=1$ for
 $j=m+1,\ldots,m+|\bsym{k}|$, we have 
\begin{align*}
 \det_{I}A
&=\frac{1}{n!}\sum_{\bsym{\s}}\sum_{\bsym{\s}_1,\ldots,\bsym{\s}_m}\prod^{m}_{j=1}\sgn(\s_j)^{\e_j}\prod^{m}_{j=1}\prod^{k_j}_{l_j=1}\sgn(\s_{j,l_j})\\
&\ \ \ \times \prod^{n}_{v=1}\a_{x_{\bsym{\s}(v)}}\bigl(x_{\bsym{\s}_{1}(v)},\ldots,x_{\bsym{\s}_{m}(v)}\bigr)\prod^{m}_{j=1}\prod^{k_j}_{l_j=1}e_{\s_{j,l_j}(v),\s_j(v)}.
\end{align*}
 Here, the sum is taken over all
 $\bsym{\s}=(\s_1,\ldots,\s_{m})\in(\fS_{n})^{m}$ and
 $\bsym{\s}_j=(\s_{j,1},\ldots,\s_{j,k_j})\in(\fS_{n})^{k_j}$ for
 $1\le j\le m$. For $\bsym{\s}=(\s_1,\ldots,\s_k)\in(\fS_{n})^k$
 and $v\in\{1,\ldots,n\}$, we set $\bsym{\s}(v):=(\s_1(v),\ldots,\s_k(v))$.
 Replacing the variables $\t_{j,l_j}=\s_{j,l_j}\s_j^{-1}$ with
 $1\le j\le m$ and $1\le l_j\le k_j$ and writing
 $\bsym{\t}_j=(\t_{j,1},\ldots,\t_{j,k_j})$, we have 
\begin{align*}
 \det_{I}A
&=\frac{1}{n!}\sum_{\bsym{\s}}\prod^{m}_{j=1}\sgn(\s_j)^{\e_j+k_j}\sum_{\bsym{\t}_1,\ldots,\bsym{\t}_m}\prod^{m}_{j=1}\prod^{k_j}_{l_j=1}\sgn(\t_{j,l_j})\\
&\ \ \ \times \prod^{n}_{v=1}\a_{x_{\bsym{\s}(v)}}\bigl(x_{\bsym{\t}_{1}\s_1(v)},\ldots,x_{\bsym{\t}_{m}\s_m(v)}\bigr)\prod^{m}_{j=1}\prod^{k_j}_{l_j=1}e_{\t_{j,l_j}\s_j(v),\s_j(v)},
\end{align*}
 where $\bsym{\t}_{j}\s_j=(\t_{j,1}\s_j,\ldots,\t_{j,k_j}\s_j)$ for
 $1\le j\le m$. Note that $e_{\t_{j,l_j}\s_j(v),\s_j(v)}=1$ holds 
 for all $1\le v\le n$ if and only if $\t_{j,l_j}=1$ since
 $x_1<x_2<\ldots<x_n$. Therefore, $\det_{I}A$ is equal to  
\begin{align*}
\frac{1}{n!}\sum_{\bsym{\s}}&\prod^{m}_{j=1}\sgn(\s_j)^{\e_j+k_j}\prod^{n}_{v=1}\a_{x_{\bsym{\s}(v)}}\bigl(x_{\bsym{\s}_1(v)},\ldots,x_{\bsym{\s}_m(v)}\bigr)\\
&=\det_{\widetilde{I}}\bigl(\a_{x_{i_1},\ldots,x_{i_{m}}}(\underbrace{x_{i_1},\ldots,x_{i_1}}_{k_1},\ldots,\underbrace{x_{i_m},\ldots,x_{i_m}}_{k_m})\bigr)_{1\le
i_1,\ldots,i_{m}\le n}.
\end{align*} 
 This completes the proof of the proposition.
\end{proof}

\begin{thm}
 Let $S=\{x_1,\ldots,x_n\}$ be a factor-closed set and
 let $F(r;n_1,\ldots,n_{k})$ be an even function $\pmod{r}$ with respect to
 the variables $n_1,\ldots,n_k$. Let $I=\{2,3,\ldots,k+1\}$ if $k$
 is even, and $\{1,2,\ldots,k+1\}$ otherwise. Then, we have
\begin{equation}
\label{for:hypdet1}
 \det_{I}\bigl(F(x_{i_1};x_{i_2},\ldots,x_{i_{k+1}})\bigr)_{1\le i_1,\ldots,i_{m+1}\le n}
=(x_1\cdots x_n)^k\prod^{n}_{v=1}\a_{x_v}(x_v,\ldots,x_v),
\end{equation}
 where $\a_{x_v}$ is the finite Fourier coefficient given by
 \eqref{for:ffcE}. 
\end{thm}
\begin{proof}
 Setting $m=1$ in Proposition~\ref{prop:hypdet}, we immediately obtain
 this theorem.
\end{proof}

\begin{remark}
 The case of $k=1$ is obtained in Theorem~$2$ in \cite{BourqueLigh1993}.
\end{remark}

\begin{cor}
\label{cor:simplestcase}
 Let $S=\{x_1,\ldots,x_n\}$ be a factor-closed set. Let
 $I=\{2,3,\ldots,m+1\}$ if $m$ is even, and $\{1,2,\ldots,m+1\}$
 otherwise. Let $\xi$ be arbitrary arithmetic function of
 $m$-variables. Then, we have
\begin{multline}
\label{for:hypdetS}
\quad \det_{I}\bigl(S^{\bsym{\g},\xi}_{\bsym{f}}(x_{i_1},\ldots,x_{i_{m+1}})\bigr)_{1\le i_1,\ldots,i_{m+1}\le n}\\
=\bigl(f_1(1)\cdots
 f_m(1)\bigr)^n\prod^{n}_{v=1}\xi(x_v,\ldots,x_v)f^{[\lcm{\g_1,\ldots,\g_{m}}]}_{m+1}(x_v).\qquad
\end{multline}
\end{cor}
\begin{proof}
 From formula \eqref{for:ffcES}, we have 
 $\a^{\bsym{\g},\xi}_{\bsym{f},x_v}(x_v,\ldots,x_v)=x_v^{-m}S^{{}^{t}\xi^{\bsym{\g}}_{x_v}}_{{}^{t}\wt{\bsym{f}}}(x_v,1,\ldots,1)$.
 Furthermore, using the identity
 $a_{\g_1}(x_v)\cdots a_{\g_m}(x_v)=a_{\lcm{\g_1,\ldots,\g_m}}(x_v)$,
 we can see that the right-hand side above is equivalent to
 $x_v^{-m}\xi(x_v,\ldots,x_v)f_1(1)\cdots f_m(1)f^{[\lcm{\g_1,\ldots,\g_m}]}_{m+1}(x_v)$.
 From formula \eqref{for:hypdet1}, we thus complete the proof.
\end{proof}

\begin{example}
 Retaining the notation in Example~\ref{ex:gcd}, we have
 $S_{\bsym{f}}=f\circ \textrm{gcd}$. Hence, by formula
 \eqref{for:hypdetS}, we have for a factor-closed set
 $S=\{x_1,\ldots,x_n\}$ 
\[
 \det_{I}\bigl((f\circ \textrm{gcd})(x_{i_1},\ldots,x_{i_{m+1}})\bigr)_{1\le i_1,\ldots,i_{m+1}\le n}
=\prod^{n}_{v=1}(f*\m)(x_v).
\]
 This is a part of the result of P. Haukkanen \cite{Haukkanen1992}.
\end{example}


\bigskip

\noindent
\textsc{Yoshinori YAMASAKI}\\
 Graduate School of Mathematics, Kyushu University.\\
 Hakozaki, Fukuoka, 812-8581 JAPAN.\\
 \texttt{yamasaki@math.kyushu-u.ac.jp}\\

\end{document}